\documentclass{amsart}
\usepackage[margin=1.5in]{geometry}

\usepackage[all]{xy}
\usepackage{tikz-cd}
\usepackage{enumerate}
\usepackage{amsfonts}
\usepackage{amssymb, amscd}
\usepackage{graphicx}
\usepackage{mathtools}
\usepackage{xcolor}
\usepackage{float}

\newcommand{\R}{\mathbb{R}}

\newcommand{\pf}{\n{\em Proof.}   }

\newcommand{\n}{\noindent}

\newcommand{\sn}{\smallskip\noindent}

\newtheorem{teo}{Theorem}[section]
\newtheorem{cor}[teo]{Corollary}

\newtheorem{lema}[teo]{Lemma}

\theoremstyle{definition}

\newtheorem{rmk}[teo]{Remark}

\title{Convex Bodies with Affinely equivalent projections and Affine Bodies of Revolution  }

\author[L. Montejano]{Luis Montejano   }
\address{IMATE, UNAM, Mexico}
\email{luis@im.unam.mx}

\date{\bf Preliminary Version}

\begin{document}
\maketitle

\begin{abstract}
In this paper, we study affine bodies of revolution. This will allow us to prove that a convex body all whose orthogonal $n$-projections are affinely
equivalent is an ellipsoid, provided $n\equiv 0,1, 2$ mod $4$, $n>1$, with the possible exemption of $n=133$. Our proof uses convex geometry and topology of compact Lie groups.
\end{abstract}

\section{Introduction}
The purpose of this paper is to prove the following Theorem. 

\begin{teo}\label{thmM}
Suppose  $B\subset \R^{N}$ is a convex body, all of whose orthogonal projections onto $n$-dimensional linear subspaces,  
for some fixed integer $n$, $1<n<N$, are affinely  equivalent.  If $n\equiv 0,1, 2$ mod $4$, then $B$ is an ellipsoid, with the possible exemption of $n=133$.
\end{teo}

The analogue for sections is the following statement which is equivalent to the real Banach conjecture \cite{B} stated en 1932.

\smallskip

  {\em Suppose  $B\subset \R^{N}$ is a convex body containing the origin in its interior, all of whose sections through $n$-dimensional linear subspaces, for some fixed integer $n$,  $1<n<N$, are affinely  equivalent.  Then $B$ is an ellipsoid},

\medskip

\noindent {\bf Real Banach Isometric Conjecture} 
{\em Suppose $V$ be a real Banach $N$-dimensional space all of whose $n$-dimensional subspaces, for some fixed integer $n$, $1< n<N$, are isometrically isomorphic to each other. Then $V$ is a Hilbert space}. 

\bigskip 

The conjecture was proved by Gromow \cite{G} in 1967,  for $n=$ even or  $N>n+1$ and by Montejano et.al. \cite{BHJM}, 
for $n\equiv 1$ mod $4$ with the possible exception of $n=133$. The history behind this conjecture can be read in \cite{So}. It is also worthwhile to see \cite{Pe} and the notes of Section 9 of \cite{MMO}. Furthermore, 
for more about bodies having affinely equivalent or congruent sections or projections,  
see Problem 3.3, Note 3.2 and Problem 7.4, Note 7.2 of Richard Gardner's Book \cite{RG}.

The reason for the strange exception $n\neq 133$ is because 133 is the dimension of the exceptional Lie group $E_7$ for which it was not possible to prove  Theorem 1.6 of \cite{BHJM}, crucial for the proof of our Theorem \ref{thmM2}.

\medskip

The proof of Theorem \ref{thmM} combines two main ingredients: convex geometry  and  topology and geometry of transformation groups.
 The first part of the article consists of using topological methods to show that, under the hypotheses of  Theorem \ref{thmM}, all orthogonal projections of $B$ are affinely equivalent to a symmetric body of revolution. Section 3 is devoted to prove the
 following characterization of ellipsoids.

\begin{teo}\label{thm:rev}
A symmetric convex body $B\subset \mathbb{C}^{n+1}$,  $n\geq 4$,  all of whose orthogonal projections onto hyperplanes are linearly equivalent to a fixed body of revolution,  is an ellipsoid.  
\end{teo}

\section{Reducing the structure group of the tangent bundle of the sphere}

During this section, we assume $B\subset \R^{n+1}$ is a convex body, all of whose orthogonal projections onto $n$-dimensional linear subspaces  are lineally equivalent, $n\equiv 0,1, 2$ mod $4$, $n>1$ and $n\not=133$.  If this is the case, we shall prove that there is a symmetric convex body of revolution $K\subset \R^{n}$ with the property that every orthogonal projection of  $B$ is linearly equivalent to $K$. 

 The link to  topology  is via a  beautiful idea that traces back to the work of Gromov \cite{G}. It consists of the following key observation.

\begin{lema}\label{lema:key} Let $B\subset \R^{n+1}$ be a symmetric convex body, all of whose orthogonal projections onto hyperplanes  are linearly equivalent to some fixed symmetric convex body $K\subset \R^n$. Let $G_K:=\{g\in SO(n)| g(K)=K\}$ be the  {\em group of symmetries} of $K$. Then the structure group of  the tangent bundle of $S^n$ can be reduced to  $G_K$. 
\end{lema}

\pf Consider the symmetric ellipsoid $E$ of minimal volume containing $B\cap \R^{n}$. For the existence of L\"owner-Johns minimal ellipsoids, see Gruber \cite{Gru}. By translation and dilatation of the principal axes of this ellipsoid, we obtain an affine isomorphism  $f:\R^{n}\to \R^{n}$ such that $f(E)$ is the unit ball of $\R^{n}$.  Let $K=f(B\cap \R^{n})$.
 Obviously,  the structure group of  the tangent bundle of $S^n$ can be reduced to  $\{g\in GL_n(\R)| g(K)=K\}$. See Section 3.1 of \cite{BHJM} for the completely analogous proof of  Lemma 1.5, for sections instead of projections, as well as a brief reminder about structure groups of differentiable manifolds and their reductions. Obviously, $K$ is affinely equivalent to every orthogonal projection of $B$. The next step is to observe that $\{g\in GL_n(\R)| g(K)=K\}$ is a subgroup of the orthogonal group 
 $O(n)$. This is so because every linear homeomorphism that fixes $K$ also fixes the unit ball. Furthermore, if the structure group of  the tangent bundle of $S^n$ can be reduced to $\{g\in O(n) | g(K)=K\}$, then it can be reduced to $G_K$.
\qed

\medskip

Lemma \ref{lema:key} can be also interpreted through the notion of  {\em  fields of convex bodies} tangent to $S^n$ as follows. See, for example, Mani \cite{Ma} and Montejano \cite{Mo1}.  For every $u\in \mathbb{S}^n$, let $u^\perp$ be the hyperplane subspace orthogonal to $u$ and $F(u)$ be the unique $n$-dimensional
ellipsoid of least volume containing the orthogonal projection $\pi_u(B)$ of $B$ in the direction $u$. The affine transformation
$\beta_u$  which maps $F(u)$ onto the $n$-dimensional Euclidean unit ball in $u^\perp$ by translating and dilating its
principal axes is a continuous function of $u$. Then all $\beta_u(\pi_u(B))$, for 
$u\in \mathbb{S}^n$ are congruent, so there is a complete turning of $\beta_{e_1}(\pi(B))$, where $\pi$ is the orthogonal projection onto 
$\R^{n}$ and $e_1=(0,\dots,0,1)$. By Lemma 2, of \cite{Mo1}, $\beta_{e_1}(\pi(B))$, is centrally symmetric, and so is any other orthogonal projection 
of $B$. This implies that

\begin{cor}\label{corA}
Let $B\subset \R^{N}$ be a convex body, all of whose orthogonal projections onto $n$-dimensional subspaces  are affinely equivalent, $1<n<N$. Then $B$ is centrally symmetric.
\end{cor}
\pf If $n+1=N$, the proof follows immediately from the above paragraph and Lemma 2 of \cite{Mo1}. The corollary follows by induction using the obvious fact that  every orthogonal projection of a convex body is centrally symmetric if and only if this body is centrally symmetric.
\qed

\medskip

As an immediate consequence  of Aleksandrov Theorem (Theorem 2.11.1 of \cite{MMO}), we have
\begin{teo}\label{thmA}
Let $B\subset \R^{N}$ be a symmetric convex body, all of whose orthogonal projections onto $n$-dimensional subspaces are congruent, $1<n<N$. Then $B$ is a ball.
\end{teo}

Let us prove now the main result of this section

\begin{teo}\label{thmM2}
Suppose  $B\subset \R^{n+1}$ is a  convex body, all of whose orthogonal projections onto hyperplanes, for some fixed integer $n$,  $1<n<N$, are affinely  equivalent.  If $n\equiv 0,1, 2$ mod $4$  and $n\not= 133$, then there exist a symmetric body of revolution $K$ such that every orthogonal projection of $B$ onto a hyperplane subspace is linearly equivalent to $K$.
\end{teo}

\pf  By the first part of this section, we now that  there exist a symmetric convex body  $K$ such that every orthogonal projection of $B$ onto a hyperplane is linearly equivalent to $K$ and such that the structure group of  $S^n$ can be reduced to  $G_K=\{g\in SO(n)\mid g(K)=K\}$. 

If $n$ is even and $n\not=6$ , Theorem 1A) of Leonard \cite{L}  implies that $G=SO(n)$ and therefore $K$ is a ball. In particular, $K$ is a body of revolution. If $n=6$, then Theorem 1A) of Leonard \cite{L} implies that $G=SU(3)$ or $U(3)$, but in both cases the action in $\mathbb{S}^5$ is transitive, which implies again that $K$ is a ball.  If $n\equiv 1$ mod $4$  and $n\not= 133$, then Theorem 1.6 of \cite{BHJM} implies that $K$ is 
a body of revolution.
\qed

\section{Affine Bodies of Revolution} 

A {\em symmetric convex body} is a compact convex subset of a finite dimensional real  vector space  with a nonempty interior, invariant under $x\mapsto -x.$  A {\em hyperplane} is a codimension 1 linear subspace. An {\em affine hyperplane} is the translation of a hyperplane by some vector.   Two sets, each a subset of a vector space,   are 
{\em linearly} (respectively, {\em affinely}) {\em equivalent} if they can be mapped to each other by a linear (respectively, affine)  isomorphism between their ambient vector spaces. Given an affine $k$-dimensional plane $\Lambda$ in $\R^n$, we denote by $\Lambda^\perp$ the corresponding 
$(n-k)$-dimensional subspace orthogonal to $\Lambda$. 

An {\em ellipsoid} is a subset of a vector space which is affinely  equivalent to the unit ball in euclidean space. 
 A convex body $K\subset\R^n$ is a {\em body of revolution} if it admits an {\em axis of revolution}, i.e.,    a 1-dimensional  line $L$ such that  each  section of $K$ by an affine hyperplane $\Delta$ orthogonal to $L$ is  an $n-1$ dimensional euclidean ball in $\Delta$, centered at $A\cap L$  (possibly empty or just a point). If $L$ is an axis of revolution of $K$ then, $L^\perp$, is the  associated  {\em hyperplane of revolution}. An {\em affine  body of revolution} is a convex body affinely  equivalent to a  body of revolution. The images, under an affine equivalence, of an  axis of revolution and its  associated hyperplane of revolution of the body of revolution are  an axis of revolution  and associated hyperplane of revolution of the affine body of revolution (not necessarily perpendicular anymore). Clearly, an ellipsoid centered at the origin is an affine symmetric body of revolution and any hyperplane serves as a hyperplane of revolution.
  
 \bigskip 
 
 The aim of this section is to prove the following theorem.

\begin{teo} \label{thmPabr}
Let $B\subset \R^{n+1}$, $n\geq 4$,  be a symmetric convex  body, all of whose   orthogonal projections are symmetric affine bodies of revolution. Then, at least one of its projections is an ellipsoid.
\end{teo} 

\bigskip

The first step is to prove that the projection of an affine body of revolution is an affine body of revolution. For that purpose, we need to prove first the following technical lemma.

\begin{lema}\label{lemsimple}
Let $\Gamma\subset \R^{n}$, $n\geq 3$, be an affine hyperplane and let $C$ be an $(n-1)$-dimensional ball contained in $\Gamma$ with center at 
$x$. Let $H$ be a hyperplane  
subspace non parallel to $H$ and let $\pi:\R^{n}\to H$ be the orthogonal projection. Then, $\pi(C)\subset H$  is an ellipsoid of revolution with axis of revolution the line  $L$, where $L=(x+ (\Gamma^\perp + H^\perp))\cap H$.  
\end{lema}
\pf  Suppose without loss of generality the center $x$ of the $(n-1)$-ball $C$ is the origin. 
Consider the linear isomorphism $f:\Gamma\to H$ given by $f(z))=\pi(z)$. Consequently, $\pi(C)=f(C)$ is an ellipsoid. Note that $f$ is the identity in the hyperplane subspace $\Gamma\cap H$ of $H$. Therefore, $\pi(C)=f(C)$ is an ellipsoid of the revolution with axis $(\Gamma^\perp + H^\perp)\cap H$, as we wished. 
\qed 
\bigskip

\begin{lema}\label{lempro}
Let $K\subset  \R^{n}$, $n\geq 3$, be an affine body of revolution with axis of revolution the line $L$ and 
let $\pi$ be an orthogonal projection along the $1$-dimensional subspace $\ell$. Then $\pi(K)$ is an affine body of revolution. 
Moreover, if $L$ is parallel to $\ell$, then $\pi(K)$ is an ellipsoid and if not, $\pi(L)$ is an axis of revolution  of  $\pi(K)$. 
\end{lema}
\pf Let us first prove the case in which $K$ is a body of revolution. We wall prove that in this case, if $L$ is parallel to $\ell$, then $\pi(K)$ is ball and if not, $\pi(K)$ is a body of revolution with axis 
$\pi(L)$. 

Suppose $L$ is not parallel to $\ell$. Let $P=\ell +L$ be the $2$-dimensional subspace generated by $\ell$ and $L$. For every $x\in L$, let $C_x=(x+L^\perp)\cap K$, then $C_x$ is either empty, the point $x$ or a ball contained in $(x+L^\perp)$ and center at $x$.  Therefore 
$$K=\cup_{\{x\in L\}} C_x, \mbox{ and }$$

\begin{equation}\label{eq:1}
\pi(K)=\cup_{\{x\in L\}} \pi(C_x).
\end{equation}

\smallskip

By Lemma \ref{lemsimple}, $\pi(C_x)$ is an ellipsoid of revolution with axis $\pi(L)$.  Consequently, by (\ref{eq:1}), $\pi(K)$ is the union of ellipsoids of revolution, all to them with the same axis $\pi(L)$. Consequently, $\pi(K)$ is a body of revolution with axis $\pi(L)$.

For the proof of the general case of the lemma we may assume, after a linear isomorphism,  that $K$ is a body of revolution and $\pi:\R^{n}\to H$ is an affine projection in the direction of $\ell$ onto the hyperplane subspace $H$ (not necessarily orthogonal to $\ell$).  Suppose $L$ is not parallel to $L$.

Consider $K\oplus\ell=\{x\in\R^{n}\mid (x+\ell)\cap K\not=\emptyset\}$. By definition $\pi(K)=(K\oplus\ell)\cap H$ and by the first part of the proof, $K_1=(K\oplus\ell)\cap\ell^\perp$ is a body of revolution with axis 
$L_1=(L + \ell)\cap\ell^\perp$. 

Let $f:\ell^\perp \to H$ be the linear map given by $f(z)=(z+\ell)\cap H$, for every $z\in \ell^\perp$. Then $f(K_1)=\pi(K)$ and $f(L_1)=\pi(L)$. Consequently $\pi(K)$ is an affine body of revolution with axis of revolution $\pi(L)$, as we wished.
\qed

\begin{lema}\label{lemun}
Let $K\subset  \R^{n}$, $n\geq 3$, be an affine body of revolution with hyperplane of revolution $H$ and 
let $\pi$ be an orthogonal projection along the line $\ell$.  Suppose that $\ell$ is parallel to $H$ and $\pi(K)$ is an ellipsoid.
Then $K$ is an ellipsoid.
\end{lema}

\pf  Since $K$ is an affine body of revolution with hyperplane of revolution $H$  and $\ell\subset H$, the shadow boundary 
$$S\partial(K,\ell)=\{x\in \mbox{bd}K\cap L'\mid L'\mbox{ is a tangent line of }K\mbox{ parallel to }\ell\}$$ 
has the following property:  there is a hyperplane $\Gamma$ such that $L\subset \Gamma$ and  $\Gamma \cap$ bd$K = S\partial(K,\ell)$.
This implies that $\pi(\Gamma \cap K)=\pi(K)$. Consequently, the section $\Gamma \cap K$ is an ellipsoid and by Lemma  2.5 of \cite{BHJM},
$K$ is an ellipsoid. 
\qed

\bigskip

From now on, let   $B\subset \mathbb{R}^{n+1}$, $n\geq4$,  be a symmetric convex  body, all of whose orthogonal projections onto  hyperplanes  are  non-elliptical, affine bodies of revolution. 

Remember that for every  line $\ell\subset \mathbb{R}^{n+1}$,  we denote by $\ell^\perp$ the hyperplane subspace of $\mathbb{R}^{n+1}$ orthogonal to $\ell$.  Furthermore, by Lemma 2.3 of \cite{BHJM}, denote by $L_\ell$ the unique  axis of revolution of the orthogonal projection of $B$ onto 
$\ell^\perp$,  by $H_\ell$ the corresponding $(n-1)$-dimensional subspace of revolution of the orthogonal projection of $B$ onto $\ell^\perp$ and 
finally denote by $N_\ell\subset \ell^\perp$  the $1$-dimensional subspace orthogonal to $H_\ell$. That is, $N_\ell=H_\ell^\perp\cap\ell^\perp$. 
Note that by the symmetry, both  $L_\ell$ and  $H_\ell$ contain the origin. 

We claim that the assignations $\ell\to L_\ell$ and $\ell\to N_\ell$ are continuos functions of $\ell$. The proof is completely analogous to the proof of Lemma 2.7 of \cite{BHJM}.

\begin{lema}\label{lemprin} Let $B\subset \mathbb{R}^{n+1}$, $n\geq 4$,  be a symmetric convex  body, all of whose  hyperplane projections are  non-elliptical, affine bodies of revolution.
Let $P$ be the plane generated by two different lines  $\ell_1$ and $\ell_2$  and suppose $N_{\ell_1}\subset \ell_2^\perp$, then 
$$\Pi_P(L_{\ell_1}))=\Pi_P(L_{\ell_2}),$$
where $\Pi_P$ is the orthogonal projection along $P$. 
\end{lema}
\pf We shall first prove that $\ell_1^\perp\cap\ell_2^\perp\cap B$ is a non-elliptical affine body of revolution. For that purpose,
first note that $\ell_1^\perp\cap\ell_2^\perp$ is orthogonal to $P$. Furthermore, $N_{\ell_1}\subset \ell_2^\perp$ implies that 
$N_{\ell_1}\subset\ell_1^\perp\cap\ell_2^\perp$ and hence by Lemma \ref{lemun}, if  $\ell_1^\perp\cap\ell_2^\perp\cap B$ is an ellipsoid, then so is 
$\ell_1^\perp\cap B$, contradicting our hypothesis.
If this is the case, by Lemma 2.3 of \cite{BHJM}, since $n-1\geq3$, we conclude that $\ell_1^\perp\cap\ell_2^\perp\cap B$ has only one axis of revolution $\Delta$.  Thus, by Lemma \ref{lempro}, considering the orthogonal projection of $\ell_1^\perp$ onto $\ell_1^\perp\cap\ell_2^\perp$ 
along $P\cap\ell_1^\perp$, we have that $\Delta=\Pi_P(L_{\ell_1})$, but on the other hand, considering the orthogonal projection of $\ell_2^\perp$ 
  onto $\ell_1^\perp\cap\ell_2^\perp$ along $P\cap\ell_2^\perp$, we have that $\Delta=\Pi_P(L_{\ell_2})$. \qed

\bigskip

Before giving the proof of Theorem \ref{thmPabr}, we need a technical topological lemma

\begin{lema}\label{lemTo}
Let $\phi,\psi:\mathbb{RP}^n\to \mathbb{RP}^n$ two continuos maps, both homotopic to the identity, $n\geq4$. 
Then $\phi^{-1}(\mathbb{RP}^{n-2})\cap \psi^{-1}(\mathbb{RP}^{n-2})\not= \emptyset$.
\end{lema}
\pf The closed set 
$\phi^{-1}(\mathbb{RP}^{n-2})$ can be thought as the carrier of $\sigma^2\in H^2(\mathbb{RP}^n)$, where $\sigma\in H^1(\mathbb{RP}^n)$ is
the generator of the cohomology ring $H^*(\mathbb{RP}^n)$. Similarly, 
$\psi^{-1}(\mathbb{RP}^{n-2})$ can be thought as  the carrier of $\sigma^2\in H^2(\mathbb{RP}^n)$. Since $n\geq4$, $\sigma^4$ is not zero and consequently 
$\phi^{-1}(\mathbb{RP}^{n-1})\cap \psi^{-1}(\mathbb{RP}^{n-2})\not= \emptyset$.
\qed

\bigskip

\noindent {\bf Proof of Theorem \ref{thmPabr}} \quad Suppose not, suppose that $B$ is a symmetric convex body all of whose orthogonal projection onto  hyperplanes sections are  non-elliptical, affine bodies of revolution.  For every  line $\ell\subset \mathbb{R}^{n+1}$, let  $H_\ell$, $L_\ell$ and $N_\ell$ as in the paragraphs before Lemma \ref{lemprin}. 

Given the line $\ell_2\subset \mathbb{R}^{n+1}$, our next purpose is to show that there exist a  line $\ell_1\subset \mathbb{R}^{n+1}$ such that

\begin{enumerate}
\item $N_{\ell_1}\subset \ell_2^\perp\cap L_{\ell_2}^\perp$, and 
\item $L_{\ell_1}\subset \ell_2^\perp\cap L_{\ell_2}^\perp.$
\end{enumerate}

For that purpose, note that the continuos  assignations $\ell\to L_\ell$ and $\ell\to N_\ell$ can be thought as continuos maps from $\mathbb{RP}^n$ into itself, moreover, since $\ell$ is orthogonal to $L_\ell$ and also orthogonal to $ N_\ell$, the assignations $\ell\to L_\ell$ and $\ell\to N_\ell$ can be thought  as continuous maps from $\mathbb{RP}^n$ into itself, which are homotopic to the identity.  By Lemma \ref{lemTo}, there exist a  line $\ell_1\subset \mathbb{R}^{n+1}$ such that
$N_{\ell_1}\subset \ell_2^\perp\cap L_{\ell_2}^\perp$ and $L_{\ell_1}\subset \ell_2^\perp\cap L_{\ell_2}^\perp.$

Note now that both $L_{\ell_1}$ and $N_{\ell_1}$ are contained in $\ell_1^\perp\cap\ell_2^\perp$.
 If $P$ is the plane 
orthogonal to $\ell_1^\perp\cap \ell_2^\perp$, 
then $N_{\ell_1}$ is orthogonal to $P$ 
and thus $P\cap \ell_1^\perp\subset H_{\ell_1}$.
By Lemma \ref{lemprin}, $\Pi_P(L_{\ell_1}))=\Pi_P(L_{\ell_2})$, where $\Pi_P$ is the orthogonal projection along $P$.

On the other hand, $L_{\ell_1}\subset\ell_1^\perp\cap\ell_2^\perp$ implies that $\Pi_P(L_{\ell_1}))=L_{\ell_1}=\Pi_P(L_{\ell_2})$. 
Finally, since 
$L_{\ell_1}$ is orthogonal to $L_{\ell_2}$ in $\ell_2^\perp$, 
then $\Pi_P(L_{\ell_2})$ is orthogonal to $L_{\ell_1}$, contradicting the fact that 
$L_{\ell_1}=\Pi_P(L_{\ell_2})$.

\qed 

\bigskip

\begin{rmk}\label{rmk:conv}
We know that the projection or section  of a convex body of revolution $B$ is again a convex body of revolution. The converse of this result, as far as we know, is an open problem. Let us state a somewhat more precise question:

\sn

{\em Let $B\subset\R^{n+1}$, $n\geq 4$, be a convex body containing the origin in its interior. Suppose every hyperplane section of B (projection onto a hyperplane)  is a body of revolution, is $B$  necessarily a body of revolution?}

\smallskip

Note that our Theorem \ref{thmPabr} points in that direction
\end{rmk}

\section{ The proof of Theorem \ref{thmM}}

The proof of Theorem \ref{thm:rev} follows immediately from Theorem \ref{thmPabr}, because by Theorem 2.12.5 of \cite{MMO}, a convex body all whose hyperplane projections are ellipsoids is an ellipsoid. 

The codimension $1$ case of 
Theorem \ref{thmM}  follows directly from Corollary \ref{corA}, Theorem \ref{thmM2} and Theorem \ref{thm:rev}. The rest of the proof follows from the fact that a convex body all whose hyperplane projections are ellipsoids is an ellipsoid. 

\bigskip 

\sn{\bf Acknowledgments.} Luis Montejano acknowledges  support  from CONACyT under 
project 166306 and  support from PAPIIT-UNAM under project IN112614.

\end{document}